\documentclass[a4paper,12pt]{article}

\usepackage[pdftex]{graphicx}
\usepackage[section] {placeins}
\usepackage{enumerate}
\usepackage[normalem]{ulem}
\usepackage{amsmath}
\usepackage{amsthm}
\usepackage{latexsym}
\usepackage{amsfonts}
\usepackage{amssymb}
\usepackage{url}

\usepackage{tikz-cd}
\usetikzlibrary{graphs}
\usetikzlibrary{arrows.meta}
\usetikzlibrary{shapes.geometric}

\usepackage{authblk}

\usepackage{latexsym,amssymb,amsfonts, amsthm, amsmath}
\usepackage[english]{babel}
\usepackage{bm}
\usepackage{mathrsfs}

\newtheorem{thm}{Theorem}[section]

\newtheorem{lem}[thm]{Lemma}
\newtheorem{cor}[thm]{Corollary}

\newtheorem{conj}[thm]{Conjecture}

\theoremstyle{definition}

\usepackage{fullpage}

\begin{document}

\title{A Comment on Dean's Construction of Prime Labelings on Ladders}

\author{Stephen J. Curran\footnote{Department of Mathematics, University of Pittsburgh at Johnstown, Johnstown, Pennsylvania, USA 15904  {\tt sjcurran@pitt.edu}} \ and   M.~A.~Ollis\footnote{Marlboro Institute for Liberal Arts \& Interdisciplinary Studies, Emerson College, Boston, Massachusetts, USA 02116  {\tt matt\_ollis@emerson.edu}}}

\date{}

\maketitle

\vspace{-10mm}
\begin{abstract}
A {\em prime labeling} on a graph of order~$m$ is an assignment of $\{ 1, 2, \ldots, m \}$ to the vertices of the graph such that each pair of adjacent vertices has coprime labels.  The {\em ladder} of order~$2n$ is the $2 \times n$ grid graph graph~$P_2 \times P_n$.  In a recent paper, Dean claimed a proof of the Prime Ladder Conjecture that every ladder has a prime labeling.   We point out a flaw in Dean's construction, showing that a stronger hypothesis is needed for it to hold.  We conjecture that this stronger hypothesis is true.  We also offer an alternative construction inspired by Dean's approach which shows that if the Even Goldbach Conjecture and a particular strengthening of Lemoine's Conjecture are true then the Prime Ladder Conjecture follows.
\end{abstract}

\section{Introduction}

Let~$G$ be a graph of order~$m$.  A {\em prime labeling} of~$G$ is an assignment of the integers $\{ 1, 2, \ldots, m \}$ to the vertices of the graph such that each pair of adjacent vertices has coprime labels.  If~$G$ has a prime labeling then it is {\em prime}.  See Gallian's dynamic survey of graph labelings~\cite{Gallian} for a comprehensive account of graph labelings in general and prime labelings in particular.

In this note we are concerned with ladders: the {\em ladder} $L_n$ of order~$2n$  is the $2 \times n$ grid graph~$P_2 \times P_n$.   The Prime Ladder Conjecture, made by Vilfred et al.~\cite{VSN02}, is that every ladder is prime.   There are several results available about the primality of ladders, unfortunately not all complete and correct.

Vilfred, Somasundaram and Nicholas show that $L_n$ is prime when~$n$, $n+1$, $n+2$ or~$2n+1$ is prime~\cite{VSN02}  and Sundaram, Ponraj and Somasundaram show that~$L_n$ is prime when~$n$ is the sum of two primes and when~$n \leq 116$~\cite{SPS07}.  Berliner~{\em et al} show that~$L_n$ is prime when~$n = (q-p)/2$ for some prime~$q$ with~$p$ either~1 or a prime smaller than~$2n+1$ \cite{B+16}.

In~\cite{Dean17}, Dean claims a full proof of the conjecture.  However, the method is flawed.  We explore this in the next section, showing that Dean's construction requires a stronger condition than the one given.  In~\cite{Schroeder19}, a proof that all cubic bipartite graphs other than $K_{3,3}$ are prime is claimed, a result that immediately implies the primality of ladders.  Unfortunately, this proof is also incorrect~\cite{SchroederErratum}.

In an unpublished work, Ghorbani and Kamali claim a proof that all ladders are prime~\cite{GhorbaniKamali}.

In this short note our contribution is two-fold.  In Section~\ref{sec:dean}  we consider what goes wrong in Dean's construction and show how the existence of ``strong canonical partitions" of integers would imply that Dean's approach is successful.  Sadly, Nathaniel Dean died in~2021 and we were unable to discuss this work with him.

In Section~\ref{sec:alt} we offer a new construction, somewhat similar to that of Dean's.   It already follows from~\cite{SPS07} that  Goldbach's Even Conjecture
 implies that~$L_n$ is prime for even~$n$.  Our construction shows that the following strengthening of Lemoine's Conjecture implies  that~$L_n$ is prime for odd~$n$.

\begin{conj}[Strengthened Lemoine's Conjecture]
For any odd integer $n\geqslant 7$, there exist primes $p$ and $q$ with $p < 2q$ such that $n=2p+q$.
\end{conj}

We have checked with GAP that the Strengthened Lemoine's Conjecture holds for $n < 10^7$.

Related work of the present authors shows that the techniques of Section~\ref{sec:alt} may be extended to the~$3 \times n$ grid graphs~$P_3 \times P_n$.  In~\cite{CO} it is shown that Goldbach's Even Conjecture and a more restrictive strengthening Lemoine's Conjecture requiring~$p<q$ are sufficient to show that these graphs are prime for all~$n$.

\section{The Flaw in Dean's Construction}\label{sec:dean}

Dean~\cite{Dean17} makes the following definition.  A {\em canonical partition} of an integer~$n$ is a sequence~$p_1, p_2, \ldots, p_m$ of odd primes such that~$\sum_{i=1}^m p_i = n$ and $p_j \geq 2\sum_{i=1}^{j-1} p_i + 3$ for~$j$ in the range~$1 < j \leq m$~\cite{Dean17}.  Dean shows that every integer~$n \geq 50$ has a canonical partition.

Dean goes on to construct a labeling for~$L_n$ given a canonical partition $p_1, \ldots, p_m$ of~$n$.  We refer the reader to~\cite{Dean17} for the full construction.  Here we note that when~$m \geq 3$ there is a potentially non-coprime adjacency unaccounted for in Dean's argument that the labeling is prime:  the elements~$\sigma_k =  2\left(\sum_{i=1}^{k-2} p_i \right)+ p_{k-1}$ and $\tau_k = 2\left(\sum_{i=1}^{k-1} p_i \right)+ p_k + 1$ are adjacent for~$k$ in the range~$3 \leq k \leq m$.

This adjacency can cause the labeling not to be prime: consider~$n=87$ with canonical partition~$3 , 11, 73$.  We obtain~$\sigma_3 =   17$ and~$\tau_3 = 102 = 6\cdot17$.   However, there are other canonical partitions for~$87$ that do work, including the one given in~\cite{Dean17}: $87 = 3 + 17 + 67$ and the partition $3,17,67$ has~$\sigma_3 = 23$ and~$\tau_3 = 108$.

Call a canonical partition $p_1, p_2, \ldots, p_m$ for~$n$ in which~$\sigma_k$ and~$\tau_k$ are coprime for all~$k$ in in the range~$3 \leq k \leq m$ {\em strong}.   The result actually proved by Dean's construction is:

\begin{thm}{\rm \cite{Dean17}}
If~$n$ has a strong canonical partition then~$L_n$ is prime.
\end{thm}

Dean conjectures that all~$n \geq 50$ have a canonical partition with at most~3 terms and confirms this up to~$n = 5,000,000$.  As a canonical partition with two terms is automatically strong, Dean's check guarantees strong canonical partitions for even~$n$ in this range.   We have checked that there is a strong canonical partition with three terms for every odd~$n$ in the same range. Therefore, we have prime labelings for~$P_2 \times P_n$ when $50 \leq n \leq 5,000,000$ via Dean's construction.

It seems likely that all~$n \geq 50$ have a strong canonical partition with at most three terms.

\section{An Alternative Construction}\label{sec:alt}

In this section we construct a prime labeling for the ladders~$L_n$, where $n = 2p+q$ for a prime~$p$ and odd prime~$q$ with $p < 2q$.
Some notation: for a label~$\ell$, let $N(\ell)$ denote the set of labels that are on adjacent vertices to the vertex labeled~$\ell$.

We start with a lemma that gives us a prime labeling for~$L_{2p}$ for prime~$p$ that has a useful additional property.

\begin{lem}\label{lem:2p}
If~$p$ is prime then there is a prime labeling for~$P_2 \times P_{2p}$ with~$1$ and~$4p$ as the labels in the final column.
\end{lem}

\begin{proof}
Here are the required labelings for~$p \in \{ 2,3,5 \}$:
{\scriptsize
\begin{center}
$\begin{array}{|c|c|c|c|}
\hline
5&4&3&8 \\
\hline
6&7&2&1\\
\hline
\end{array}$
\hspace{7mm}
$\begin{array}{|c|c|c|c|c|c|}
\hline
7 & 2 & 3 & 10 & 11 & 12 \\
\hline
6 & 5 & 4 & 9 & 8 & 1 \\
\hline
\end{array}$
\hspace{7mm}
$\begin{array}{|c|c|c|c|c|c|c|c|c|c|}
\hline
15 & 2 & 3 & 4 & 17 & 14 & 5 & 18 & 19 & 20 \\
\hline
16 & 7 & 8 & 9 & 10 & 11 & 12 & 13 & 6 & 1 \\
\hline
\end{array}$
\end{center} }

So assume~$p \geq 7$.   Consider the labeling~$\mathcal{S}$ given by
$$f(u_{i,j}) = \begin{cases}
  (i-1)p + j & \text{when } 1 \leq j \leq p, \\
  (3-i)p + j & \text{when } p+1 \leq j \leq 2p.\\
  \end{cases}$$
There are two pairs of neighboring non-coprime labels in~$\mathcal{S}$: in column~$p$ we have the entries~$p$ and~$2p$ and in column~$2p$ we have the entries~$3p$ and~$4p$.

Construct a new labeling by switching labels 1 and $3p$ and switching labels 4 and $2p$.   Any non-coprimality must involve one of these labels.

We have~$N(1) = \{ 3p-1, 4p \}$, $N(4) = \{ p, 2p-1, 2p+1 \}$, $N(2p) = \{ 3, 5, p+4 \}$ and $N(3p) = \{ 2,  p+1 \}$.   The only potential non-coprimality is between~$3p$ and~$p+1$, which share a factor of~3 when~$p \equiv 2 \pmod{3}$.  In this instance, switch the labels~$p$ and~$3p$.  Then $N(p) = \{ 2, p+1 \}$ and $N(3p) = \{ 4, p-1,  3p+1 \}$ and no non-coprime neighbours remain.
\end{proof}

The following array is the prime labeling for~$L_{22}$ given by Lemma~\ref{lem:2p}:
{\scriptsize
$$\begin{array}{|c|c|c|c|c|c|c|c|c|c|c|c|c|c|c|c|c|c|c|c|c|c|}
\hline
11 &2 & 3 & 22 & 5 & 6 & 7 & 8 & 9 & 10 & 33 & 34 & 35 & 36 & 37 & 38 & 39 & 40 & 41 & 42 & 43 & 44 \\
\hline
12 & 13 & 14 & 15 & 16 & 17 & 18 & 19 & 20 & 21 & 4 & 23 & 24 & 25 & 26 & 27 & 28 & 29 & 30 & 31 & 32 & 1 \\
\hline
\end{array}$$ }

\begin{thm}\label{th:2p+q}
If~$n = 2p+q$ where~$p$ is prime and~$q$ is an odd prime, then~$L_n$ is prime.
\end{thm}

\begin{proof}
Let $\mathcal{S}_1$ be the prime labeling on $L_{2p}$ given in Lemma~\ref{lem:2p}.
Extend the labeling $\mathcal{S}_1$ to the labeling
$\mathcal{S}_2$ on $L_{2p+q}$
by letting
$f(u_{i,j}) = 2p+(i-1)q+j$ for all
$1 \leqslant i \leqslant 2$ and
$2p+1 \leqslant j \leqslant 2p+q$.
Let $k=\lfloor 4p/q \rfloor$ and let column $j^*$ be the column with the labels $\{ (k+1)q, (k+2)q \}$.
The internal labels in column $j^*$ are the only non-coprime
labels in~$\mathcal{S}_2$. We show that we can swap one of the labels in column~$j^*$ with another to produce a new labeling~$\mathcal{S}_3$ that is prime.

First suppose that either~$p < q$ or~$\frac12 p < q < \frac23 p$.  If~$2p < q$ then~$2q$ is in column~$j^*$, if $p < q < 2p$ then~$4q$ is in column~$j^*$ and if $\frac12 p < q < \frac23 p$ then~$8q$ is in column~$j^*$.  Therefore, as all adjacent labels in~$\mathcal{S}_2$ are of opposite parity,  it is sufficient to find a power of~2 that is not adjacent to a multiple of~$q$ to swap with~$2p$, $4p$ or $8p$ and let~$\mathcal{S}_3$ be the labeling resulting from this swap.

If $p=2$ then $N(8) = \{ 1,3,9 \}$ which we can use when~$q \neq 3$.  If $q = 3$ then there is no suitable power of~2 as~3 is in each of $N(2)$, $N(4)$ and~$N(8)$.  However, $N(12) = \{ 1,9,13 \}$ and $N(14) = \{11, 13 \}$ and so we may swap~$12$ and~14 to obtain a prime labeling.

If~$p = 3$ then $N(4) = \{3,5,9\}$, which is suitable when $q \not\in \{3,5\}$.  We have $q \neq 3$.  If $q=5$ then  $N(8) = \{1, 9, 11 \}$ so we may use~8.

If $p=5$ then $N(8) = \{ 3,7,9 \}$ which is suitable for~$q \not\in \{ 3,7 \}$.  If $q = 3$ then no suitable power of~2 is available.  Instead, swap~$7q = 21$, which is in column~$j^*$ and has neighbors $N(21) = \{20, 22, 24 \}$, with~23, which has neighbors $N(23) = \{ 22, 26 \}$.  If $q = 7$, then~$N(4) = \{ 3, 9, 17 \}$.

If $p \geq 7$ then $N(2) = \{ 3, p, p+2 \}$ or $\{ 3, p+2, 3p \}$.  As $q \neq 3$ and $q \neq p$, this is suitable for all~$q$ except $q = p+2$.   If $q = p+2$ then $p>7$ and we have $N(8) = \{ 7,9, p+8 \}$.

Now suppose that $\frac23 p < q \leq p$.    In this range, if $\frac45 p  < q$ then the entries in column~$j^*$ are $\{ 5q, 6q \}$ and if $\frac45 p  >  q$  then the entries in column~$j^*$ are $\{ 6q, 7q \}$.   It is sufficient to find a label of the form~$2^a3^b$ with $a,b > 0$ that is not adjacent to a multiple of~$q$ to swap with~$6p$ and let~$\mathcal{S}_3$ be the labeling resulting from this swap.

We cannot have $p=2$.  If $p=3$ then~$q=3$ and we can use $N(6) = \{ 5,7 \}$.  If $p=5$ then $q \leq 5$ and we can use $N(6) = \{1, 13, 19 \}$.

If $p=7$ then $q \leq 7$ and $N(12) = \{5,11,13 \}$, leaving only $q=5$ to consider.  There is no label of the form~$2^a3^b$ that we can use in this case.  Instead, we switch the other entry in column~$j^*$, which is $5p = 35$, with $p=7$.  We have $N(7) = \{4,6,22 \}$ and $N(35) = \{ 30, 34, 36 \}$.

If $p \geq 11$ then $q > \frac23 p \geq 7$ and  $N(6) = \{ 3,5, p+6 \}$.   If $q \mid p+6$ then $p+6$ must be a non-trivial odd multiple of~$q$.  Hence $p+6 \geq 3q > 2p$ and so $p<6$, which it is not and hence $\gcd(q,p+6) = 1$.
\end{proof}

The following array is the prime labeling for $P_2 \times P_{21}$ given by Theorem~\ref{th:2p+q} using $p=5$ and $q=11$.
{\scriptsize
$$\begin{array}{|c|c|c|c|c|c|c|c|c|c|c|c|c|c|c|c|c|c|c|c|c|}
\hline
15 & 2 & 3 & 4 & 17 & 14 & 5 & 18 & 19 & 20 &
21 & 8 & 23 & 24 & 25 & 26 & 27 & 28 & 29 & 30 & 31  \\
\hline
16 & 7 & 22 & 9 & 10 & 11 & 12 & 13 & 6 & 1 &
32 & 33 & 34 & 35 & 36 & 37 & 38 & 39 & 40 & 41 & 42  \\
\hline
\end{array}
$$
}

Upon noting that $L_1$, $L_3$ and $L_5$ are known to be prime~\cite{VSN02}, we reach our claimed result.

\begin{cor}
If the Strengthened Lemoine's Conjecture holds then $P_2 \times P_n$ is prime for all odd~$n$.
\end{cor}

\end{document}